\documentclass[a4paper,11pt]{article}
\usepackage{amsmath,amscd,amsfonts,amssymb,epsf,latexsym}
\usepackage[all]{xy}

\makeatletter

\long\def\@makefnt#1{\parindent 1em\noindent
            \hb@xt@1.8em{\hss\@textsuperscript{}}#1}
\long\def\@ftntext#1{\insert\footins{%
    \reset@font\footnotesize
    \interlinepenalty\interfootnotelinepenalty
    \splittopskip\footnotesep
    \splitmaxdepth \dp\strutbox \floatingpenalty \@MM
    \hsize\columnwidth \@parboxrestore
    \color@begingroup
      \@makefnt{%
        \rule\z@\footnotesep\ignorespaces#1\@finalstrut\strutbox}%
    \color@endgroup}}%
\def\subjclass#1{%
  \@ftntext{2010 {\itshape Mathematics Subject Classification.}\enspace #1.}}
\def\keywords#1{%
  \@ftntext{{\itshape Key words and phrases.}\enspace #1.}}
\makeatother

\def\moins{\raise 1pt\hbox{{$\scriptstyle -$}}}
\def\plus{\raise 1pt\hbox{{$\scriptstyle +$}} }

\newtheorem{theorem}{Theorem}
\newtheorem{proposition}[theorem]{Proposition}
\newtheorem{lemma}[theorem]{Lemma}
\newtheorem{corollary}[theorem]{Corollary}
\newtheorem{remark}[theorem]{Remark}

\newtheorem{note}[theorem]{Note}

\def\proof{\noindent{\bf Proof.\ }}

\def\qed{~\hbox{$\Box$}}

\def\rank{\mathop{\rm rank}}
\def\dim{\mathop{\rm dim}}

\def\Hom{\mathop{\rm Hom}}

\def\Ker{\mathop{\rm Ker}}

\def\cO{{\mathcal O}}

\def\cX{{\mathcal X}}

\def\Z{{\mathbb Z}}

\def\C{{\mathbb C}}

\def\Q{{\mathbb Q}}

\begin{document}

\title{\bf Diagonals of flag bundles}

\author{Shizuo Kaji\thanks{Research partially supported by KAKENHI, Grant-in-Aid for Young Scientists (B) 26800043 and JSPS Postdoctoral Fellowships for Research abroad. }\\
\small Department of Mathematical Sciences, Faculty of Science, Yamaguchi University\\
\small 1677-1, Yoshida, Yamaguchi 753-8512, Japan\\
\small skaji@yamaguchi-u.ac.jp \and
Piotr Pragacz\thanks{Supported by National Science Center (NCN)
grant no. 2014/13/B/ST1/00133}\\
\small Institute of Mathematics, Polish Academy of Sciences\\
\small \'Sniadeckich 8, 00-656 Warszawa, Poland\\
\small P.Pragacz@impan.pl}



\date{}

\maketitle

%\tableofcontents

\begin{abstract}
We express the diagonals of projective, Grassmann and, more generally, flag bundles of type $(A)$
using the zero schemes of some vector bundle sections, and do the same for their single point subschemes.
We discuss diagonal and point properties of these flag bundles. We study when the complex manifolds $G/B$ 
for other groups have the point and diagonal properties. We discuss explicit formulas for the classes 
of diagonals of the varieties $G/B$.
\end{abstract}

\section{Introduction}
For a map of varieties $\pi: F\to X$, it is useful to study the diagonal
in the fibre square $F\times_X F$.
The classes of such diagonals for fibre bundles in topology and smooth proper morphisms in algebraic geometry
were investigated by Graham \cite{G}, Fulton and the second author \cite{P1}, \cite[Appendix G]{FP}. As explained in \cite[Section 5]{P1},
knowing such a class, one can compute the class of a subscheme of $F$. For an overview of applications, see \cite[Chapter 7]{FP}.

In the present paper, we shall rather study the
diagonals in the Cartesian squares $F\times F$ of the total spaces of flag bundles $\pi: F\to X$. Suitable resolutions of the structure
sheaves of the diagonals over the structure sheaves of the Cartesian squares of some homogeneous spaces were used by Kapranov in \cite{Ka} to give 
descriptions of their derived categories. Many schemes can be realized as degeneracy loci of vector bundle homomorphisms. It turns out that 
to understand degeneracy loci, it is useful to study diagonals of flag bundles (cf. \cite{P1}, \cite{DP}, \cite{F2}, \cite{PR}, \cite{FP}).

In \cite{PSP}, Pati, Srinivas and the second author investigated which varieties $X$ have
the following ``diagonal property'' ($D$): there exists a vector bundle of rank $\dim(X)$
on $X\times X$ with a section whose zero scheme is the diagonal. 
If $X$ has $(D)$, then it is nonsingular.

Also, the following ``weak point property'' ($P$)
was investigated: for some point $x\in X$, there exists a vector bundle of rank $\dim(X)$ on $X$ with a section 
whose zero scheme is $x$. If a variety $X$ has $(D)$, then it has $(P)$ for any $x\in X$. Any nonsingular curve has $(D)$. The product of varieties
having $(D)$, has $(D)$. In \cite[Section 3]{PSP}, we gave several detailed results on surfaces with $(D)$.
In particular, it was shown ({\it loc.cit.}, Proposition 4) that a ruled surface (cf., e.g., \cite[Chap. V, Sect. 2]{Ha}) has $(D)$, i.e., the projectivization of any rank 2 vector bundle on a nonsingular curve has $(D)$. This result was one of the starting points of the present paper.

It was shown by Fulton \cite{F}
and the second author \cite{P} that the flag varieties of the form $SL_n/P$ over any field have $(D)$. 
Samuelsson and Sepp\"anen \cite{SS} gave recently an application of the diagonal property of flag varieties to global complex analysis.
The interest to $(D)$ for flag varieties was related to the theory of Schubert polynomials
of Lascoux-Sch\"utzenberger \cite{LS}. These authors defined on the polynomial ring in several variables a scalar product 
\cite[Chapter 10]{L}, 
for which the (single) Schubert polynomials and their ``duals'' form adjoint bases ({\it loc.cit.}, Corollary 10.2.4). The reproducing kernel for this scalar product is equal to the top (double)
Schubert polynomial ({\it loc.cit.}, Section 10.2), which, in turn, is equal to the top Chern class of the vector bundle, realizing $(D)$ for the variety
$SL_n/B$ of complete flags -- a result of Fulton \cite{F}. 

\smallskip

A natural question emerged then: do the flag varieties for other groups have the diagonal property?
(see \cite[p. 1235]{PSP}, \cite[Conjecture 8.2]{P} and \cite[Proposition 12]{PSP}). In the present paper, we answer this question almost completely for the manifolds $G/B$, where $G$ is a simple, simply connected, complex algebraic group (see Section \ref{full}). We show that for $G$ of type $(B_i) (i\ge 3), (D_i) (i\ge 4), (G_2), (F_4)$ and $(E_i) (i=6,7,8)$, the flag manifold $G/B$ has not the diagonal property (see Theorem \ref{thm:point-property-for-flag}).
The main tools are the Atiyah-Hirzebruch homomorphism
and the Borel characteristic homomorphism. 
In Section \ref{form}, we recall several explicit formulas for the classes of diagonals of $G/B$; we do not use, however, these formulas in
the proof of Theorem \ref{thm:point-property-for-flag}.

\smallskip

Apart from the question about generalized flag manifolds, the present paper arose from our attempts to understand relations 
between the diagonals of the base spaces and those of the total spaces of flag bundles.

We shall also study the following variant of $(D)$. 
Given a section $s$ of a vector bundle, we write $Z(s)$ for its zero scheme.
We say that $X$ has property $(D')$, if there
exist two vector bundles $A$ and $B$ on $X\times X$ such that $\rank(A) + \rank(B)=\dim(X)$, a section $s$ of $A$ and 
a section $t$ on $Z(s)$ of 
the restriction $B_{Z(s)}$ of $B$ to $Z(s)$ such that $Z(t)$ is the diagonal of $X$. 
Thus for $A=(0)$, we recover $(D)$. Note that $(D')$ is a slight weakening of $(D)$ as the key property
that the rank of the bundle is $\dim(X)$ holds also for $(D')$. A variety $X$ with $(D')$ is nonsingular because 
its cotangent sheaf is locally free: it is isomorphic to the restriction of $A^{\vee}\oplus B^{\vee}$ to the diagonal of $X$. 

The main results of the paper are Theorems \ref{D'}, \ref{pf}, \ref{p'f} and \ref{thm:point-property-for-flag}.

The paper is organized as follows.
In Section \ref{vb}, we discuss properties of flag bundles and construct a certain vector bundle
on the Cartesian square of a flag bundle. Using this vector bundle, we prove in Theorem \ref{D'} 
that if the base space of a flag bundle has $(D)$, then its total space has $(D')$.
In Section \ref{top}, we discuss the analogs of this result for topological properties $(D_r)$ and $(D_c)$ from
\cite[Section 6]{PSP}. The topological situation is easier than that in algebraic geometry: if $X$ has
$(D_r)$ or $(D_c)$ then the corresponding flag bundles also have these properties (see Theorem {\ref{drE}).

In Section \ref{point} we show that if a quasiprojective base of a flag bundle has $(P)$, then its total space has $(P)$
(see Theorem \ref{pf}). For more general schemes, we prove in Theorem \ref{p'f} that if the base space of a flag bundle has $(P)$, 
then its total space has $(P')$, a property analogous to $(D')$. In Section \ref{top}, we discuss the analogs for topological properties $(P_r)$ 
and $(P_c)$ from \cite[Section 6]{PSP}. 

In Section \ref{full}, we investigate which complex manifolds $G/B$ for other groups $G$ have $(P_c)$, see
Theorem \ref{thm:point-property-for-flag} and Corollary \ref{red}. An absence of $(P_c)$ implies for many of them
the absence of $(D_c)$.

In Appendix, we discuss explicit formulas for the classes of diagonals of the varieties $G/B$ due to Fulton, the second author and Ratajski,
Graham, and De Concini, adding one for the type $(G_2)$ (see Lemma \ref{ng2}). We also disprove an integrality conjecture from \cite{G}
(see Remark \ref{cg}).

\section{Flag bundles of type $(A_{n-1})$}\label{vb}

Let $E$ be a vector bundle of rank $n$ on a variety $X$ over a field.
Fix an increasing sequence of integers
$$
d_{\bullet}: 0<d_1<d_2<\ldots<d_{k-1}<d_k=n\,.
$$
By a ${d_{\bullet}}$-flag, we mean an increasing sequence of subbundles of $E$
$$
E_1\subset E_2 \subset \cdots \subset E_{k-1}\subset E_k=E
$$
such that $\rank(E_i)=d_i$ for $i=1,\ldots,k$. Let 
$$
\pi: Fl_{d_{\bullet}}(E)\to X
$$ 
be the flag bundle parametrizing all ${d_{\bullet}}$-flags. 
For example, the sequence
$$
d_1=d<d_2=n
$$
gives rise to the Grassmann bundle $G_d(E)$, parametrizing subbundles of rank $d$ of $E$ (see \cite[B.5.7]{Fit}).
For $d=1$, we get the projectivization of $E$: $P(E)=G_1(E)$ (see \cite[Section 7]{Ha} and \cite[B.5.5]{Fit}).

It is well-known that
\begin{equation}\label{dimG}
\dim(G_d(E))=\dim(X)+d(n-d)\,.
\end{equation}

On $Fl_{d_{\bullet}}(E)$, there exists 
the following tautological sequence of vector bundles:
\begin{equation}\label{tauto}
{S}_1\hookrightarrow {S}_2\hookrightarrow \cdots \hookrightarrow {S}_{k-1}
\hookrightarrow {S}_k={\pi}^*(E) \stackrel{q_{1}}\twoheadrightarrow {Q}_1 \stackrel{q_{2}}
\twoheadrightarrow {Q}_2 \stackrel{q_{3}}\twoheadrightarrow \cdots \stackrel{q_{k}}\twoheadrightarrow
{Q}_k=0\,,
\end{equation}
where $\rank({S}_i)=d_i$ for $i=1,\dots,k$, and
${Q}_i$ is the quotient of $\pi^*{E}$ by ${S}_i$\,,
so that $\rank({Q}_i)=n-d_i$.
On $G_d(E)$, the tautological sequence (with $S=S_1$, $Q=Q_1$)
\begin{equation}\label{exa}
0\to S\to \pi^*E \to Q\to 0\,,
\end{equation}
where $\rank(S)=d$, is a short exact sequence.

\smallskip

Regarding $Fl_{d_{\bullet}}(E)\to X$ as a tower of Grassmann bundles
\begin{equation}\label{tower}
Fl_{d_{\bullet}}(E)=G_{n-d_{k-1}}(Q_{k-1})\to \cdots \to G_{d_3-d_2}(Q_2) \to  G_{d_2-d_1}(Q_1)\to G_{d_1}(E)\to X\,,
\end{equation}
and using (\ref{dimG}), we see that with $d_0=0$, we have
\begin{equation}\label{dim}
\dim (Fl_{d_{\bullet}}(E))= \dim(X)+\sum_{i=1}^{k-1}(d_i-d_{i-1})(n-d_i)\,.
\end{equation}

Write $F=Fl_{d_{\bullet}}(E)$.
Let $F_1=F_2=F$, and denote by
$$
p_i:F_1\times F_2 \to F_i
$$
the two projections. We shall now construct a certain vector bundle of rank $\dim(F)-\dim(X)$ on $F_1\times F_2$.
If $k=2$, using the notation of (\ref{exa}), we define the following vector bundle
\begin{equation}\label{H2}
H= {\Hom}(p_1^*S,p_2^*Q)=(p_1^*S)^\vee \otimes p_2^*Q\,.
\end{equation}
Suppose now that $k\ge 3$. Using the notation of (\ref{tauto}),
consider the following homomorphism of vector bundles on $F_1\times F_2$:
$$
\varphi: \bigoplus_{i=1}^{k-1} 
{\Hom}(p_1^*{S}_i,p_2^*{Q}_i)
\to \bigoplus_{i=1}^{k-2} 
{\Hom}(p_1^*{S}_i,p_2^*{Q}_{i+1})\,,
$$
defined by
\begin{equation}\label{phi}
\varphi(\sum_{i=1}^{k-1} h_i)=\sum_{i=1}^{k-2} \bigl(h_{i+1}|p_1^*{S}_i - p_2^*(q_{i+1})\circ h_i\bigr)\,,
\end{equation}
where $h_i\in 
{\Hom}(p_1^*{S}_i,p_2^*{Q}_i)$.

\begin{lemma} The homomorphism $\varphi$ is surjective.
\end{lemma}
\proof
Let us fix $i=1,\ldots, k-2$. 
Let $h\in 
{\Hom}(p_1^*{S}_i,p_2^*{Q}_{i+1})$. By subtracting from $h$ a suitable 
homomorphism from ${\Hom}(p_1^*{S}_{i+1},p_2^*{Q}_{i+1})$
restricted to the subbundle $p_1^*{S}_i$ of $p_1^*S_{i+1}$, we get a homomorphism from $p_1^*S_i$ to $p_2^*Q_{i+1}$, 
which factorizes through $p_2^*{Q}_i$.
But such a homorphism belongs to $\varphi(
{\Hom}(p_1^*{S}_i,p_2^*{Q}_i))$. The assertion follows.
\qed

\smallskip

Define the following vector bundle on $F_1\times F_2$:
\begin{equation}\label{H}
H=\Ker (\varphi)\,.
\end{equation}
Using Lemma 1, we obtain (with $d_0=0$)
\begin{equation}\label{rH}
\rank (H)=\sum_{i=1}^{k-1} d_i(n-d_i) - \sum_{i=1}^{k-2} d_i(n-d_{i+1})=\sum_{i=1}^{k-1}(d_i-d_{i-1})(n-d_i)\,.
\end{equation}

\begin{note} \rm The present section is an expanded version of \cite[pp. 107-8]{P}. The bundle $H$, defined in (\ref{H}), is modeled on the bundle $K$ from \cite[(7.6)]{F}.
\end{note}

\begin{remark} \rm Let $X$ be a point. It is shown in \cite[Section 7]{F} that for 
$$
d_{\bullet}=0<1<2<\ldots <n-1<n\,,
$$ 
the top Chern class of $H$ is the top
double Schubert polynomial taken on first Chern classes of the tautological quotient bundles on the two copies of
complete flag varieties. It will follow from Section \ref{diag} that it is actually the class of the diagonal of a complete flag variety.
\end{remark}

\section{Diagonal properties}\label{diag}

We adopt the set-up from the previous section, and state the following result.

\begin{theorem}\label{D'} \ If $X$ has $(D)$, then for any vector bundle $E$ and any $d_{\bullet}$, $Fl_{d_{\bullet}}(E)$ 
has $(D')$.
\end{theorem}
\proof 
Let $G$ be a vector bundle of rank $\dim(X)$ on $X\times X$ with a section whose zero scheme $Z(s)$ is the diagonal $\Delta_X$ of $X$.
Fix $d_{\bullet}$, and follow the notation from Section \ref{vb}. 
Let
$$
G'=(\pi_1\times\pi_2)^*(G)
$$
be a bundle on $F_1\times F_2$ together with a section $s'=(\pi_1\times\pi_2)^*(s)$. 
Consider 
$$
Z:=Z(s')=(\pi_1\times \pi_2)^{-1}(\Delta_X) \subset F_1\times F_2\,.
$$
Let $r_1, r_2: X\times X\to X$ be the two projections. The following two vector bundles on $\Delta_X$ 
are equal:
\begin{equation}\label{qEX}
(r_1^*E)_{\Delta_X}=(r_2^*E)_{\Delta_X}\,.
\end{equation}
Since 
$$
(\pi_1\times\pi_2)^* r_i^*E=p_i^*(E_{F_i})
$$
for $i=1,2$, we obtain from (\ref{qEX}) that the following two vector bundles on $Z$ are equal:
\begin{equation}\label{id}
\bigl(p_1^*E_{F_1}\bigr)_{Z}=\bigl(p_2^*E_{F_2}\bigr)_{Z}\,.
\end{equation}
Thanks to (\ref{id}), we get, for any $i=1,\ldots,k-1$, the following homomorphism:
\begin{equation}\label{ti}
h_i: \bigl(p_1^*S_i\bigr)_Z \to \bigl(p_1^*E_{F_1}\bigr)_{Z}=\bigl(p_2^*E_{F_2}\bigr)_{Z} \to \bigl(p_2^*Q_i\bigr)_Z
\end{equation}
of vector bundles on $Z$. Here, the subbundle $S_i\hookrightarrow E_F$
and the quotient bundle $E_F \twoheadrightarrow Q_i$ are from (\ref{tauto}).
The family of homorphisms $\{h_i\}$ gives rise to the section 
$$
h=\sum h_i\in \Gamma \bigl(Z,\oplus_{i=1}^{k-1} \Hom(p_1^*S_i, p_2^*Q_i)_Z\bigr)\,.
$$

Suppose $k\ge 3$. It follows from (\ref{ti}) that we have on $Z$
$$ 
h_{i+1}|p_1^*S_i=p_2^*(q_{i+1})\circ h_i
$$
for $i=1,\ldots,k-2$. Indeed, since $h_i$ and $h_{i+1}$ factorize through the bundle (\ref{id}), the two homomorphisms
$$
h_{i+1}|p_1^*S_i \ , \ \ p_2^*(q_{i+1})\circ h_i \ : \ \ \bigl(p_1^*S_i\bigr)_Z \to \bigl(p_2^*Q_{i+1}\bigr)_Z
$$
are equal.
Invoking (\ref{phi}), we see that 
$$
\varphi \circ h=0\,,
$$ 
so $h$ induces
a section $t$ of the bundle $H_{Z}$, where $H$ is the vector bundle on $F_1\times F_2$ from (\ref{H}) and
(\ref{H2}).
By (\ref{dim}) and (\ref{rH}), we have 
$$
\rank(G')+\rank(H)=\dim(F)\,.
$$

We claim that the section $t$ of the bundle $H_Z$
vanishes precisely (scheme theoretically) on the diagonal $\Delta_F \subset F_1\times F_2$. It vanishes on $\Delta_F$ since
the tautological sequence of vector bundles on $G_{d_i}(E)$ is a complex for any $i=1,\ldots,k-1$ (cf. (\ref{exa})).

Having defined the sections $s', t$ globally, it is sufficient to check the converse assertion $Z(t)\subset \Delta_F$ locally, where
$F_1\times F_2$ is the product of the Cartesian square of the base space times the Cartesian square of the flag variety 
$Fl_{d_{\bullet}}(E_x)=: F_x$, where $x\in X$. In other words, this boils down to check the assertion over the point $x$, i.e. on $F_x\times F_x$.
For the case of complete flags, see \cite[p. 402]{F}. For any $d_{\bullet}$, let $f\in Z$ with $\pi_1(x)=\pi_2(x)=x$, so we may 
regard $f$ as a point 
$$
f=(L_1 \subset \cdots \subset L_{k-1}\subset L_k=E_x \ , \ M_1\subset \cdots \subset M_{k-1}\subset M_k=E_x)
$$ 
in $F_x\times F_x$.
Let $F_{x,1}=F_{x,2}=F_x$. For any $i=1,\ldots,k-1$, the restriction of $h_i$ (see \ref{ti}) to $F_{x,1}\times F_{x,2}$ is
\begin{equation}\label{tV}
p_1^*S_i \to p_1^*V_{F_{x,1}}=V_{F_{x,1}\times F_{x,2}}= p_2^*V_{F_{x,2}} \to p_2^*Q_i\,,
\end{equation}
where we write $V$ for $E_x$, $S_i$ and $Q_i$ are the restrictions to $F_x$ of the tautological bundles on $F$, and $p_1, p_2$ are the two projections from $F_{x,1}\times F_{x,2}$
to the factors.
At the point $f=((L_i),(M_i))$, (\ref{tV}) becomes the map
$$
L_i\hookrightarrow V \twoheadrightarrow V/M_i\,,
$$
whose vanishing implies $L_i=M_i$. This holds for any $i=1,\ldots,k-1$. We have proved that set-theoretically $Z(t)=\Delta_F$.
It is not hard to verify that this equality holds scheme-theoretically.
The assertion of the theorem follows.
\qed

\smallskip

We record the following simple fact.

\begin{lemma}\label{triv} \ Let $E$ be a vector bundle on a variety $X$. Let $r_1, r_2: X\times X \to X$ be the two projections.
Suppose that the following two vector bundles on $X\times X$ are equal:
$$
r_1^*E=r_2^*E\,.
$$
Then $E$ is a trivial bundle.
\end{lemma}
\proof
Fix a point $x\in X$. By the assumption, we have
$$
(r_1^*E)_{X\times \{x\}}=(r_2^*E)_{X\times \{x\}}\,.
$$
Via the identification $X\times \{x\}\simeq X$, the LHS is the bundle $E\to X$. The RHS is the trivial bundle $(E_x)_X$. The assertion follows.
\qed

\begin{remark}\label{analiza} \rm 
Let us speculate a bit about {\it this} proof of Theorem \ref{D'}.
To convert it to that of $(D)$, we must extend the section $t$ to the whole $F_1\times F_2$. This can be done
only if $p_1^*(E_{F_1})=p_2^*(E_{F_2})$; so, by virtue of Lemma \ref{triv}, only if the bundle $E_F$ is trivial.
\end{remark}

\section{Point properties}\label{point}

We first record the following result.

\begin{proposition}\label{pg} Suppose that $X$ is a quasiprojective variety with $(P)$.
Then for any vector bundle $E$ on $X$ and any $1\le d\le n-1$, the Grassmann bundle $G_d(E)$ has $(P)$.
\end{proposition}
\proof
By the assumption, for a certain point $x\in X$, there exists a vector bundle $G$ of rank $\dim(X)$ and 
$s\in \Gamma(X,G)$ such that $Z(s)=x$.

Suppose that $\rank(E)=n$. We realize $X$ as an open subset in a projective variety $X'$. By \cite[Proposition 2]{BS},
there exists a coherent sheaf $E'$ on $X'$ whose restriction to $X$ is $E$. Let $L$ be the restriction of $\cO_{X'}(1)$ to $X$.

We claim that there exists an integer $m$ such that $E\otimes L^{\otimes m}$ has $n$ global sections which
are independent at $x$.
Indeed, this follows (by restriction from $X'$ to $X$) from \cite[Th\'eor\`eme 2(a), p. 259]{Se} which asserts 
that there exists an integer $m$ such that the $\cO_{x,X'}$-module $E'(m)_x$ is generated by the elements of $\Gamma(X',E'(m))$. 
Choose any $d$ sections out of these $n$ global sections of $E\otimes L^{\otimes m}$.
Using a canonical isomorphism 
$$
G_d(E\otimes L^{\otimes m})\simeq G_d(E)\,,
$$
and the assumption on $E$, we can reduce to the situation when $E$ has $d$ sections $\{s_i\}$ which are independent at $x$. 

Set $F=G_d(E)$, and denote by $\pi: F\to X$ the projection. Let $Q$ be the tautological quotient rank $n-d$
bundle on $F$. Consider the following rank $\dim(F)=\dim(X) + d(n-d))$ vector bundle $H$ on $F$:
$$
H=\pi^*G \oplus Q^{\oplus d}\,.
$$
We define the following section $t\in \Gamma(F,H)$. On the first summand of $H$, we take the pullback via $\pi^*$ of 
the section $s: X \to G$. On the last $d$ summands, we take the following sections: we compose
the sections 
$$
\pi^*(s_i): F \to \pi^*E
$$ 
with the canonical surjection $\pi^*E \twoheadrightarrow Q$. We have
\begin{equation}\label{zt}
Z(t)=Z(\pi^*(s))\cap Z(\oplus \pi^*(s_i))\,.
\end{equation}
But
$$
Z(\pi^*(s))=\pi^{-1}(x)=G_d(E_x)\,,
$$
so that (\ref{zt}) as a point of $G_d(E_x)$ corresponds to the $d$-dimensional vector subspace of $E_x$ spanned by $(s_i)_x$.
We get that $Z(t)$ is a single point in $G_d(E)$. Hence $G_d(E)$ has $(P)$.
\qed

\begin{remark} \rm The projective bundle 
$$
P\bigl(\cO(a_1)\oplus \cdots \oplus \cO(a_r)\bigr)\to P^n\,,
$$
$a_i\in \Z$, is a toric variety (cf. \cite{O}). Thus the proposition
gives some support to the conjecture (cf. \cite[p. 115]{P}) that a nonsingular toric variety has $(P)$ (and
perhaps even $(D)$ -- which is known in the surface case \cite{PSP}).
\end{remark}

It is well-known (cf., e.g., \cite[p. 142]{Gr}) that the projective, Grassmann and flag bundles
on quasiprojective varieties are quasiprojective. Realizing a flag bundle as a tower (\ref{tower})
of Grassmann bundles, and using an easy induction, we infer from the proposition the following
result.

\begin{theorem}\label{pf} Suppose that $X$ is a quasiprojective variety with $(P)$. Then for any
vector bundle $E$ on $X$ and any $d_{\bullet}$, the flag bundle $Fl_{d_{\bullet}}(E)$ has $(P)$.
\end{theorem}

\begin{remark} \rm 
Sometimes, one
studies the following ``strong point property'' of a variety $X$: for {\it any} $x\in X$, there exists a bundle
on $X$ of rank $\dim(X)$ with a section whose zero scheme is $x$. Granting this property for a quasiprojective variety $X$, 
the above reasoning shows 
that $F=Fl_{d_{\bullet}}(E)$ has $(P)$ for any point $f\in F$. Indeed, given $f\in F$, put $x=\pi(f)$, and argue as above.
Thus $F$ also has the strong point property.
\end{remark}

For {\it any} scheme, we still have a result explained in the following theorem. We need a definition.
We say that $X$ has
property $(P')$, if for some $x\in X$, there
exist two vector bundles $A$ and $B$ on $X$ such that $\rank(A) + \rank(B)=\dim(X)$, a section $s$ of $A$ and 
a section $t$ of $B_{Z(s)}$ such that $Z(t)$ is $x$. 
If $X$ has $(D')$, then for any $x\in X$, $(P')$ holds by restricting the data giving $(D')$ to $X\times \{x\}$.

\begin{theorem}\label{p'f} \ If $X$ has $(P)$, then $Fl_{d_{\bullet}}(E)$ has $(P')$ for any $d_{\bullet}$.
\end{theorem}
\proof
Fix $d_{\bullet}$, and write $F=Fl_{d_{\bullet}}(E)$. Let $\pi: F\to X$ be the projection.
Suppose that for the fixed point $x\in X$, there exists a vector bundle $G$ of rank $\dim(X)$ on $X$ with a section $s$ whose zero scheme is $x$. 
We shall show that $F$ has $(P')$ for any point $f\in \pi^{-1}(x)$. Let $G'=\pi^*(G)$ and $s'=\pi^*(s)$. 
Consider $W=Z(s')\subset F$. In other words $W=\pi^{-1}(x)$.

Using the vector bundle $H$ from (\ref{H}) and (\ref{H2}), we define the following vector bundle:
$$
H'=H_{F \times \{f\}}
$$ 
on $F \simeq F\times \{f\}$.
Note that $\rank(G')+\rank(H')=\dim(F)$.
Invoke $Z=(\pi\times \pi)^{-1}(\Delta_X)\subset F\times F$ from the proof of Theorem \ref{D'}. 
In this proof, we constructed the section $t$ of the bundle $H_Z\to Z$ whose zero scheme
is the diagonal of $F$.
We have $W\simeq W\times \{f\} \subset Z$.
The restriction to $W\times \{f\}$ of the section $t$, gives rise to a section, denoted $t'$, of the bundle
$H'_{W\times \{f\}}\to W\times \{f\}$. 
We claim that $Z(t')=f$. The section $t'$ vanishes at $f$ because $(f,f)$ belongs to the diagonal.

It is sufficient to check the converse assertion locally. Let $g\in Z(t')$.
Since $\pi(g)=\pi(f)=x$, we may regard
$$
f=(L_1\subset \cdots \subset L_{k-1}\subset L_k=E_x) \ \ \ \hbox{and} \ \ \ g=(M_1\subset \cdots \subset M_{k-1}\subset M_k=E_x)
$$
as points in $Fl_{d_{\bullet}}(E_x)=F_x$.
Write $V=E_x$. For $i=1,\ldots,k-1$, we consider (\ref{tV}) restricted to $F_x \times \{f\}$:
\begin{equation}\label{hres}
p_1^*S_i \to p_1^*V_{F_x} = V_{F_x\times \{f\}}= p_2^*V_f \to p_2^*(Q_i)_f \,,
\end{equation}
where $p_1: F_x\times \{f\} \to F_x$, $p_2: F_x\times \{f\} \to f$ are the two projections, and $S_i$ (resp. $Q_i$) are the restrictions of
the tautological bundles from $F$ to $F_x$ (resp. $f$).
Restricted to the point $g$, (\ref{hres}) becomes the map
$$
M_i\hookrightarrow V \twoheadrightarrow V/L_i\,,
$$
whose vanishing implies $M_i=L_i$.
This holds for every $i=1,\ldots,k-1$. We have proved that $g=f$, i.e., $Z(t')=f$, and hence $F$ has $(P')$ for any $f\in \pi^{-1}(x)$.
\qed

\begin{remark} \rm
Granting the strong point property for $X$, the above reasoning shows 
that $F$ has $(P')$ for any point $f\in F$. Indeed, given $f\in F$, put $x=\pi(f)$, and argue as above.
\end{remark}

\section{Topological properties}\label{top}

We now pass to topology. We first recall some definitions from \cite[Section 6]{PSP}.
Let $X$ be a (smooth) compact connected oriented manifold, and $\Delta$ be the diagonal submanifold
of $X\times X$. We say that $X$ has property ($D_r$) if there exists a smooth real vector bundle of rank $\dim(X)$
on $X\times X$ with a smooth section $s$ which 
is transverse to the zero section of the bundle and whose zero locus is $\Delta$.
If $\dim_{\mathbb{R}} X=2m$ and the above vector bundle is a complex vector bundle
of complex rank $m$, then we say that $X$ has property ($D_c$). If $X$ has
($D_c$), then it is almost complex ({\it loc.cit.}, p. 1259). For a complex manifold, we have the following relation
between the diagonal properties: \ $(D) \Rightarrow (D_c) \Rightarrow (D_r)$.

\begin{remark} \rm In \cite[Section 6]{PSP}, the diagonal property ($D_o$) is also studied (one requires that the
bundle involved in the definition of ($D_r$) is orientable). It is proved there that a real projective
space of odd dimension does not have ($D_o$).
\end{remark}

Let $E$ be a smooth real vector bundle on $X$. For any $d_\bullet$ like in Section \ref{vb}, there is an associated flag bundle 
$\pi: Fl_{d_{\bullet}}^{\mathbb R}(E)\to X$
parametrizing $d_{\bullet}$-flags of real subbundles of $E$. It is endowed with the tautological sequence (\ref{tauto})
of real bundles. Similarly,
if $E$ is a smooth complex vector bundle on $X$, then there is an associated flag bundle $\pi: Fl_{d_{\bullet}}^{\mathbb C}(E)\to X$
parametrizing $d_{\bullet}$-flags of complex subbundles of $E$, endowed with the tautological sequence (\ref{tauto})
of complex bundles.

\begin{theorem}\label{drE} \ (i) If $X$ has ($D_r$) and $E\to X$ is a smooth real bundle, then $Fl_{d_{\bullet}}^{\mathbb R}(E)$
has ($D_r$).

\noindent
(ii) \ If $X$ has ($D_c$) and $E\to X$ is a smooth complex bundle, then $Fl_{d_{\bullet}}^{\mathbb C}(E)$ has ($D_c$).
\end{theorem}
Proof.
Both cases of the theorem can be proved by the construction using the tautological bundles from the proof
of Theorem \ref{D'}. Using the proof of this theorem and its notation, we have
$$
\Delta_F \subset Z \subset F_1 \times F_2\,,
$$
and we have the section $s'$ of $G'$ and the section $t$ of $H_Z$.
By a partition of unity argument (cf. \cite[Lemma 1.4.1]{A}),
$t$ can be extended to a global section of $H$.
Then, $s' \oplus t$ is a global section of $G'\oplus H$ which vanishes exactly on $\Delta_F$
(compare with Remark \ref{analiza}).
\qed

We say, following \cite[Section 6]{PSP}, that $X$ as above has property ($P_r$) if there exists a smooth real vector bundle of rank $\dim(X)$
on $X$ with a smooth section $s$ which is transverse to the zero section of the bundle and whose zero locus is a point.
If $\dim_{\mathbb{R}} X=2m$ and the above vector bundle is a complex vector bundle
of complex rank $m$, then we say that $X$ has property ($P_c$). 

\begin{remark}\label{no} \rm It was shown in \cite[Remark 6]{PSP} that if a bundle $E$
of rank $\dim(X)$ on a (connected) manifold $X$ has $e(E)=\pm 1$ (resp. $c_m(E)=\pm 1$), then we can use this bundle to realize $(P_r)$ (resp. $(P_c)$).
\end{remark}

\begin{remark} \rm
By the argument from Theorem \ref{p'f},
we see that if $X$ has $(P_r)$ (resp. $(P_c)$), then $Fl_{d_{\bullet}}^{\mathbb R}(E)$ (resp. $Fl_{d_{\bullet}}^{\mathbb C}(E)$)
has  $(P_r)$ (resp. $(P_c)$). The same holds for the corresponding strong point properties.
\end{remark}

\section{Manifolds $G/B$ for other groups}\label{full}

A general reference for group-theoretic notions used in this section
is \cite{Hu}. 
Let $G$ be a simple, simply connected algebraic group over $\C$,
$B$ its Borel subgroup, and $T$ a maximal torus contained in $B$.
Denote by $G/B$ the generalized flag manifold.

In this section, we work in topological category and all vector bundles are complex.
Suppose that the complex dimension of $G/B$ is $m$. We shall study 
when $G/B$ has $(P_c)$, i.e. when there exists a vector bundle $E$
of complex rank $m$ on $G/B$ such that $c_m(E)$ is the class of a point in $H^{2m}(G/B;\Z)$
(cf. Remark \ref{no}).

Our main result in this section is
\begin{theorem}\label{thm:point-property-for-flag}
For $G$ of type $(B_i) (i\ge 3), (D_i) (i\ge 4), (G_2), (F_4)$ and $(E_i) (i=6,7,8)$,
the flag manifold $G/B$ has not $(P_c)$, and consequently it has not the diagonal property $(D_c)$.
\end{theorem}

To prove the theorem, we need several results.
Let $\cX(T)$ be the group of characters of $T$ and let $K(G/B)$ be the Grothendieck group of $G/B$ (cf. \cite[Section 2.1]{A}).
Consider the {\it Atiyah-Hirzebruch homomorphism} (see \cite[Definition 3.17(a)]{KK}):
$$
\beta_1: S(\cX(T)) \to K(G/B)
$$
such that for $\lambda\in \cX(T)$, $e^\lambda\mapsto \ \hbox{class of} \ L_\lambda=G\times_B \C_\lambda$,
a line bundle on $G/B$.
Here, we regard the $T$-representation $\C_\lambda$ as a $B$-representation
by letting the nilradical of $B$ act trivially. Then, we extend this definition multiplicatively to the entire symmetric algebra $S(\cX(T))$.

We record (see {\cite[Theorem 4.6]{KK}} and the references therein):
\begin{theorem}
The homomorphism $\beta_1$ is surjective.
\end{theorem}
Since in $S(\cX(T))$ any element is a $\Z$-linear combination of monomials $e^{\lambda_1}\cdots e^{\lambda_k}$, where $\lambda_i\in \cX(T)$, and

\begin{equation}
\begin{split}
\beta_1(e^{\lambda_1}\cdots e^{\lambda_k})=\beta_1(e^{\lambda_1})& \cdots \beta_1(e^{\lambda_k})\\
=&[L_{\lambda_1}]\cdots [L_{\lambda_k}]= [L_{\lambda_1}\otimes \cdots \otimes L_{\lambda_k}]=[L_{{\lambda_1}+\cdots +{\lambda_k}}]\,,
\end{split}
\end{equation}
the theorem implies the following

\begin{corollary}\label{cor:splitting}
In $K(G/B)$, the class of any vector bundle is a $\Z$-linear combination of the classes of line bundles $L_\mu$ for some $\mu \in \cX(T)$.
\end{corollary}

\begin{remark} \rm
It is shown by Kumar in \cite[Corollary 2.12]{K} that for $G/B$ the present $K$-group is the same as the algebraic geometric $K$-group 
discussed in \cite[Section 15.1]{Fit}.
\end{remark}

Recall now the following {\em Borel characteristic homomorphism}:
\[
c: S(\cX(T)) \to H^*(G/B;\Z)
\]
such that for $\lambda \in \cX(T)$, $e^\lambda\mapsto c_1(L_\lambda)$. Then, we extend this definition multiplicatively to all $S(\cX(T))$
(see  \cite{Bo1} and \cite{Dem} for more details).

It follows from Corollary \ref{cor:splitting} that

\begin{corollary}\label{cor:image-of-c^*}
The Chern classes of any vector bundle on $G/B$ are in the image of $c$.
\end{corollary}

The smallest positive integer $t_G$ such that 
$$
t_G \cdot \text{(class of a point)}
$$
is in the image of $c$ is called the {\em torsion index} of $G$.
We record (see {\cite[Proposition 4.2]{Bo}, \cite[Proposition 7]{Dem} and also \cite{T}):
\begin{theorem}\label{thm:torsion-index}
We have $t_G=1$ if and only if $G$ is of type $(A_i)$ or $(C_i)$.
\end{theorem}

Combining Corollary \ref{cor:image-of-c^*} and Theorem \ref{thm:torsion-index},
the assertion of Theorem \ref{thm:point-property-for-flag} follows.\qed

Let $G$ be a complex reductive group.
Recall that by replacing $G$ with its universal covering group, the flag variety $G/B$
can be regarded as a product of flag varieties associated to simple, simply-connected groups of type 
$(A_i)$, $(B_i)$ ($i \ge 3$), $(C_i)$ ($i \ge 2$), $(D_i)$ ($i \ge 4$), $(G_2)$, $(F_4)$ and $(E_i)$ ($i = 6, 7, 8$).

\begin{corollary}\label{red} \rm Let $G$ be a complex reductive group containing either type $(B_i)$($i \ge 3$), $(D_i)$ ($i \ge 4$), $(G_2)$, $(F_4)$ 
or $(E_i)$ ($i = 6, 7, 8)$ as a factor.
Then, its flag variety $G/B$ has not $(P_c)$, and consequently it has not $(D_c)$.
\end{corollary}
\proof
The class of a point in $G_1/B_1 \times G_2/B_2$ is the product of the classes of points in $G_1/B_1$ and $G_2/B_2$,
and 
$$
K(G_1/B_1\times G_2/B_2)\simeq K(G_1/B_1)\otimes K(G_2/B_2)
$$ 
by the K\"unneth theorem. Therefore, our argument for the proof of Theorem 18 applies straightforwardly.\qed

\smallskip

We pass now to type $(C_n)$. We have an identification of $F=Sp(2n,\mathbb C)/B$ with the space of complete isotropic flags 
$$
V_1\subset V_2\subset \cdots \subset V_n \subset \C^{2n}\,.
$$
Let 
$$
(0)=S_0\subset S_1\subset S_2\subset \cdots \subset S_n\subset \C^{2n}_F
$$
be the tautological flag on $F$. We set for $i=1,\ldots,n$, $L_i=S_i/S_{i-1}$. Let $x_i=c_1(L_i)$. 
Then, we have
$$
H^*(F;\Z)=\dfrac{\Z[x_1,x_2,\ldots,x_n]}{(e'_i(x))},
$$
where $e'_i(x)$ is the $i$-th elementary symmetric polynomial in $x_1^2,x_2^2,\ldots, x_n^2$.
Now, the class of a point is $x_1 x_2^3 \cdots x_n^{2n-1}$, which is the top Chern class of
the bundle
$$
L_1 \oplus L_2{^{\oplus 3}} \oplus \cdots \oplus L_n{^{\oplus 2n-1}}\,.
$$
This shows that the following result holds.
\begin{proposition}\label{p} The flag manifold $Sp(2n,\C)/B$ has $(P_c)$.
\end{proposition}

\begin{remark} \rm In the paper \cite{KP}, we shall investigate the diagonal and point
properties of the spaces $G/P$.
\end{remark}

\section{Appendix: Explicit formulas for the classes of diagonals}\label{form}

Several authors worked out algebraic expressions for the classes of diagonals of the spaces $G/B$ in the cohomology rings $H^*(G/B\times G/B;\Z)$
(or equivalently in the Chow rings $A^*(G/B\times G/B)$). 
Let us mention them and relevant references: Fulton \cite{F, F1}, the second author and Ratajski \cite{PR}, Graham \cite{G}, 
De Concini \cite{DeC} (see also \cite{FP}). 

Let $BG$ and $BB$ denote the classifying spaces of $G$ and $B$.
Consider the sequence 
$$
G/B \times G/B \to BB \times_{BG} BB \to BB\times BB\,,
$$
which yields the following sequence of homomorphisms of their cohomology rings:
$$
S(\cX(T))\otimes S(\cX(T))\to H^*(BB\times_{BG} BB;\Z)\to H^*(G/B\times G/B;\Z)\,.
$$
The first map is the Borel characteristic homomorphism and it is surjective after tensoring with $\Q$. 
Thus we shall realize the representatives of the classes of diagonals of $G/B$ in
\begin{equation}\label{SS}
S(\cX(T))_{\Q}\otimes_{\Q} S(\cX(T))_{\Q}=\Q[x_1,\ldots,x_n; y_1,\ldots,y_n]\,,
\end{equation}
where $x_i\in \cX(T)$ and $y_i\in \cX(T)$ are coordinates on $T\times T$\footnote{Our convention is that $x_i$ and $y_i$ represent the same coordinate
on $T$; moreover, on $T\times T$, $x_i$ corresponds to the first factor and $y_i$ to the second.}.

In \cite[Theorem 1.1]{G}, the author established a criterion for an element
in (\ref{SS}) to represent the class of the diagonal $\Delta$ of $G/B$. Among other methods, we mention Gysin maps (\cite{P1, PR, G, FP}) and equivariant cohomology (\cite{DeC}, cf. also the end of the introduction to \cite{G}).

\smallskip

For type $(A_{n-1})$, $[\Delta]$ is represented by 
$$
\prod_{1\le i<j\le n}(x_i-y_j)\,,
$$ 
the top (double) Schubert polynomial of Lascoux-Sch\"utzenberger (see \cite{F}). (Note that
our convention for numbering the variables is different from that in \cite{F}.)

\smallskip

For type $(C_n)$, $[\Delta]$ is represented by 
$$
\prod_{i<j}(x_i-y_j)\cdot W(x,y)\,,
$$
where denoting by $e_i(x)$ the $i$th elementary symmetric polynomial in $x$,
$$
W(x,y)=|e_{n+1+j-2i}(x)+e_{n+1+j-2i}(y)|_{1\le i,j \le n}\,
$$
(see \cite{F1}), or (see \cite{PR}):
$$
W(x,y)=\sum_{I\subset (n,\ldots,1)} \widetilde Q_I(x)\cdot \widetilde Q_{(n,\ldots,1) \smallsetminus I}(y)\,,
$$ 
where the sum is over strict partitions $I$, and $\widetilde Q_I(x)$ is defined in \cite[Section 4]{PR}. In \cite{DeC},
the following expression representing $[\Delta]$ was given:
\begin{equation}\label{concini}
\prod_{i<j}(x_i^2-y_j^2)\prod_i(x_i+y_i)\,.
\end{equation}

For type $(B_n)$ we have analogous expressions (differing by powers of 2) (cf. \cite{F1, F2, DeC, PR}). This is also the
case of the formulas for type $(D_n)$ from \cite{F1, PR}. 
The author of \cite{DeC} stated that $[\Delta]$ is represented by $$
\prod_{i=1}^{n-1}W_i(x,y)\,,
$$
where
$$
W_i(x,y)=\frac{1}{2}\bigl((1+\frac{y_i}{x_i})\prod_{j>i}(x_i^2-y_j^2)+(-1)^{n-i}(x_i\cdots x_n-y_i\cdots y_n)\frac{y_{i+1}\cdots y_n}{x_i}\bigr)\,. 
$$ 
We can regard this expression in the following way.
First, we multiply $W_i(x,y)$ by $x_i$. Then the first summand is the one for type $(B_n)$, which is almost good but of degree by one greater.
So we add a zero class so that it makes the sum divisible by $x_i$. This gives rise to the second summand.
Notice that the only term in the first summand which is not divisible by $x_i$ is $y_i y_{i+1}^2 \cdots y_n^2$.
Hence we have to add another term involving the $y_i$'s to cancel out the term in cohomology. This gives the above expression.

\smallskip

For type $(G_2)$, in \cite{G}, the following expression was given for a representative of $[\Delta]$:
$$
-\frac{27}{2}(x_1-y_2)(x_1-y_3)(x_2-y_3)(x_1 x_2 x_3 + y_1 y_2 y_3).
$$
Note that here $S(\cX(T))=\Z[x_1,x_2,x_3]/(x_1+x_2+x_3)$.

\begin{remark}\label{cg} \rm
In \cite[p. 483]{G}, the author conjectured that 
$$
\frac{1}{2}(x_1x_2x_3+y_1y_2y_3)
$$ 
is integral. A Schubert calculus computation gives that this class is a $\Q$-linear combination of Schubert classes
\begin{equation}\label{dispr}
-\frac{2}{9} \sigma_{s_2 s_1 s_2} + \hbox{lower terms}\footnote{Here, we use the identification
$
H^*(BB \times_{BG} BB) = H^*_T(G/B)
$
with the bi-grading on the RHS given by the degree in the coefficient $H^*(BT)$ and the degree of the Schubert basis.
So ``lower terms'' mean a $H^*(BT)$-linear combination of Schubert classes corresponding to shorter Weyl group elements.}
\end{equation}
The expression (\ref{dispr}) disproves the conjecture.
\end{remark}

We add another expression for the type $(G_2)$. We now identify
$$
S(\cX(T))_\Q \otimes_\Q S(\cX(T))_\Q=\Q[a_1,a_2;b_1,b_2]\,,
$$
where $a_1$ and $b_1$ (resp. $a_2$ and $b_2$) the two copies of the simple short root (resp. long root).
\begin{proposition}\label{ng2} The expression
\begin{equation}\label{3g2}
\frac12 (a_1+b_1)(a_1-(2b_1+b_2))(a_1+(2b_1+b_2))(a_1-(b_1+b_2))(a_1+(b_1+b_2))(a_2-(3b_1+b_2))
\end{equation}
represents $[\Delta]$.
\end{proposition}
\proof
We use Graham's criterion \cite[Theorem 1.1]{G}.
The Weyl group of $G_2$ is the dihedral group of order $12$.
The orbit of $a_1$ under $W$ is 
$$
\{\pm a_1, \pm (2a_1+a_2), \pm(a_1+a_2) \}\,.
$$
This accounts for the first five factors.
The stabilizers of $a_1$ are the identity and $s_{3a_1+2a_2}$, 
and the latter takes $a_2$ to $3a_1+a_2$. This accounts for the last factor.
If we evaluate the polynomial (\ref{3g2}) at $a_1=b_1, a_2=b_2$, we obtain the product of all the positive roots.
By Grahams's criterion, the expression (\ref{3g2}) represents $[\Delta]$.
\qed

\noindent
A similar expression was stated in \cite{DeC}. This representative 
$$
\frac12 (a_1+b_1)(a_1^2-(2b_1+b_2)^2)(a_1^2-(b_1+b_2)^2)(a_2-(3b_1+b_2))
$$
has property that the global coefficient equals the inverse of the torsion index, which is the best possible.

%\begin{remark} \rm
%In the papers  \cite{LP1, LP2}, Lascoux and the second author compared the results of \cite{F1, F2} with those of \cite{PR}. 
%\end{remark}

\vskip10pt
\noindent
{\bf Acknowledgments} 

\vskip5pt

We thank Adrian Langer for helpful discussions.

A part of the present article was written during the stay of the second author at IMPA in Rio de Janeiro, in May 2015.
He thanks this institute, and especially Eduardo Esteves, for hospitality.

\end{document}